\documentclass[a4paper,8pt]{article}
\usepackage{amsmath}
\usepackage{amsthm,amsfonts,amssymb}
\usepackage[dvips]{graphicx}
\usepackage[english]{babel}
\usepackage{color}
\usepackage{graphics}
\usepackage{graphpap}

\providecommand{\U}[1]{\protect\rule{.1in}{.1in}}

\providecommand{\U}[1]{\protect\rule{.1in}{.1in}}
\newtheorem{theorem}{Theorem}

\makeatletter
\def\blfootnote{\xdef\@thefnmark{}\@footnotetext}
\makeatother

\begin{document}

\title{Analytic solution for preferential attachment probabilities scheme}
\author{Andrea Monsellato}
\maketitle

\begin{abstract}
\noindent
Preferential attachment probabilities scheme appear in the context of scale free random graphs \cite{Barabási},\cite{Dorogovtsev}.
In this work we present preferential attachment probabilities scheme as a sequence of dependent Bernoulli random variables and we give an explicit expression of the transition probabilities.
\end{abstract}


\section{Introduction}

\noindent
Preferential attachment probabilities schemes arise in the context of scale free random graphs, following the description used in \cite{Dorogovtsev}:\\

\noindent
\emph{The preferential linking arises in this simple model not because of some special rule including a function of vertex
degree as in \cite{Barabási},\cite{Dorogovtsev 1} but quite naturally. Indeed, in the model that we consider here, the probability that a
vertex has a randomly chosen edge attached to it is equal to the ratio of the degree $k$ of the vertex and the total number of edges, $2t-1$.}\\

\noindent
we want to mimic this dynamic using a sequence of dependent Bernoulli random variables.\\

\noindent
Let $\{Y_m(n)\}_{n>0}$ a sequence of r.v. such that $Y_m(n)\sim Ber\left(\frac{\sum_{i=m}^{n-1}Y_m(i)}{2n-1}\right)$, with $Y_m(1)=1$. The sequence $\{Y_m(n)\}_{n>0}$, represent the trails to attach a new arcs at node $m$ when a new
node $n$ is generated and jointed to the network, obviously $n>m$. The degree of the node $m$ will be $X_m(n)=\sum_{i=m}^n Y_m(i)$.\\

\newpage
\noindent
We recover transition probabilities using elementary probability arguments, then

\begin{small}
\begin{align*}
P(X_m(n)=k)=&P(X_m(n)=k|X_m(n-1)=k-1)P(X_m(n-1)=k-1)+\\&
+P(X_m(n)=k|X_m(n-1)=k)P(X_m(n-1)=k)=\\&
\\&
=P(X_m(n-1)+Y_m(n)=k|X_m(n-1)=k-1)P(X_m(n-1)=k-1)+\\&
+P(X_m(n-1)+Y_m(n)=k|X_m(n-1)=k)P(X_m(n-1)=k)=\\&
\\&
=P(Y_m(n)=1|X_m(n-1)=k-1)P(X_m(n-1)=k-1)+\\&
+P(Y_m(n)=0|X_m(n-1)=k)P(X_m(n-1)=k)=\\&
\\&
=\frac{k-1}{2n-1}P(X_m(n-1)=k-1)+\frac{2n-1-k}{2n-1}P(X_m(n-1)=k)
\end{align*}
\end{small}

\noindent
Now set $p^m_{n,k}=P(X_m(n)=k)$, we have that

\begin{align}\label{General Master equation}
p^m_{n,k}=p^m_{n-1,k-1}\frac{k-1}{2n-1}+p^m_{n-1,k}\frac{2n-1-k}{2n-1},\quad 0\leq m < n,\quad 0\leq k\leq n-m
\end{align}

\section{Recovering explicit expression of transition probabilities}

\noindent
Choosing $m=0$ in (\ref{General Master equation}), i.e. we consider the probabilities transition equations for the first node, we have 

\begin{align}\label{PA Master Equation p(n,k)}
\left\{
\begin{array}{ll}
p_{1,1}=1,\quad p_{n,k}=0,\quad \forall k>n\\
p_{n,k}=\frac{k-1}{2n-1}p_{n-1,k-1}+\frac{2n-1-k}{2n-1}p_{n-1,k},\quad \forall k\leq n
\end{array}
\right.
\end{align}

\noindent
Observe that (\ref{General Master equation}) is time invariant for all nodes scaling properly the quantity $n-m$.\\

\noindent
For (\ref{PA Master Equation p(n,k)}) holds that:

\begin{small}
\begin{align*}
p_{n,1}=\frac{2(n-1)}{2n-1}p_{n-1,1}\quad \forall n\geq 2
\end{align*}

\noindent
i.e.

\begin{align}\label{PA Master Equation Sol p(n,1)}
p_{n,1}=\frac{2^{n-1}(n-1)!}{\prod_{i=1}^{n}(2i-1)}\quad \forall n\geq 1
\end{align}

\noindent
also

\begin{align*}
p_{n,n}=\frac{n-1}{2n-1}p_{n-1,n-1}\quad \forall n\geq 2
\end{align*}

\noindent
i.e.

\begin{align*}
p_{n,n}=\frac{(n-1)!}{\prod_{i=1}^{n}(2i-1)}\quad \forall n\geq 1
\end{align*}
\end{small}

\noindent
Let

\begin{align}\label{PA Master Equation Sol p(n,n)}
a_{n,k}=\frac{2^{n-k}\prod_{i=1}^{n}(2i-1)}{(n-1)!}p_{n,k}
\end{align}

\noindent
from (\ref{PA Master Equation p(n,k)}),(\ref{PA Master Equation Sol p(n,1)}) and (\ref{PA Master Equation Sol p(n,n)}) we have that

\begin{align}\label{PA Master Equation a(n,k)}
\left\{
\begin{array}{ll}
a_{n,1}=2^{2(n-1)},\quad a_{n,n}=1,\quad a_{n,k}=0 \quad \forall k>n\\
a_{n,k}=\frac{1}{n-1}\left[(k-1)a_{n-1,k-1}+2(2n-1-k)a_{n-1,k}\right]\quad \forall k=1,...,n-1\quad n>1
\end{array}
\right.
\end{align}

\noindent
Now we give an analytic solution of (\ref{PA Master Equation a(n,k)}).\\ 

\begin{theorem}
Let $\{a_{n,k}\}$ as in (\ref{PA Master Equation a(n,k)}), then

\begin{align}\label{PA Master Equation Sol a(n,k)}
\left\{
\begin{array}{ll}
a_{n,n}=1\\
a_{n,k}=2^{2(n-k)}\left[1+(k-1)\sum_{j=0}^{n-k-1}\frac{2^{-2(j+1)}(k+2j)!}{(j+1)!(k+j)!}\right]\quad \forall k=1,...,n-1\quad n>1
\end{array}
\right.
\end{align}

\noindent
\textbf{Proof}\\

\end{theorem}

\noindent
If $k=1$ obvious.\\

\noindent
If $k=n-1$ and $\forall n\geq 2$, then (\ref{PA Master Equation a(n,k)}) will be

\begin{align}\label{PA Master Equation a(n,n-1)}
a_{n,n-1}=\frac{1}{n-1}\left[(n-2)a_{n-1,n-2}+2na_{n-1,n-1}\right]
\end{align}

\noindent
now substituting (\ref{PA Master Equation Sol a(n,k)}) in (\ref{PA Master Equation a(n,n-1)}) we have that

\begin{align*}
&2^{2}\left[1+(n-2)\sum_{j=0}^{0}\frac{2^{-2(j+1)}(n-1+2j)!}{(j+1)!(n-1+j)!}\right]=\\&
=\frac{n-2}{n-1}2^{2(n-1-n+2)}\left[1+(n-3)\sum_{j=0}^{n-1-n+2-1}\frac{2^{-2(j+1)}(n-2+2j)!}{(j+1)!(n-2+j)!}\right]+\frac{2n}{n-1}\\&
\\&
2^2+(n-2)=\frac{n-2}{n-1}2^{2}\left[1+(n-3)2^{-2}\right]+\frac{2n}{n-1}\\&
\\&
2^2(n-1)-2n=2n-4\\&
\\&
-4=-4
\end{align*}

\newpage
\noindent
Let $k=2,...,n-2$ substituting we have that

\begin{align}{\label{PA Master Equation Sol a(n,k) k=2,...,n-2}}
\notag &(n-1)\left[1+(k-1)\sum_{j=0}^{n-k-1}\frac{2^{-2(j+1)}(k+2j)!}{(j+1)!(k+j)!}\right]=(k-1)\left[1+(k-2)\sum_{j=0}^{n-k-1}\frac{2^{-2(j+1)}(k-1+2j)!}{(j+1)!(k-1+j)!}\right]\\&
\notag\\&
\notag (n-1)\left[1+(k-1)\sum_{j=0}^{n-k-1}\frac{2^{-2(j+1)}(k+2j)!}{(j+1)!(k+j)!}\right]=\frac{1}{2}(2n-1-k)\left[1+(k-1)\sum_{j=0}^{n-k-2}\frac{2^{-2(j+1)}(k+2j)!}{(j+1)!(k+j)!}\right]
\notag\\&
(n-1)\left[1+(k-1)\sum_{j=0}^{n-k-1}\frac{2^{-2(j+1)}(k+2j)!}{(j+1)!(k+j)!}\right]=\left[(n-1)-\frac{1}{2}(k-1)\right]\left[1+(k-1)\sum_{j=0}^{n-k-2}\frac{2^{-2(j+1)}(k+2j)!}{(j+1)!(k+j)!}\right]
\end{align}

\noindent
we recover that (\ref{PA Master Equation Sol a(n,k) k=2,...,n-2}) is equivalent to

\begin{align}{\label{PA Master Equation Sol Final a(n,k) k=2,...,n-2}}
&\sum_{j=0}^{n-k-2}\frac{2^{-2(j+1)}}{(j+1)!}\left[2(k-2)\frac{(k-1+2j)!}{(k-1+j)!}-(k-1)\frac{(k+2j)!}{(k+j)!}\right]
=2^{-2(n-k-1)}\frac{(2n-k-3)!}{(n-k-1)!(n-2)!}-1
\end{align}

\noindent
Now let $n=k+2+r$, $r\geq 0$ the equality (\ref{PA Master Equation Sol Final a(n,k) k=2,...,n-2}) become

\begin{align}{\label{PA Master Equation Sol Final v1  a(n,k) k=2,...,n-2}}
&\sum_{j=0}^{r}\frac{2^{-2(j+1)}}{(j+1)!}\left[2(k-2)\frac{(k-1+2j)!}{(k-1+j)!}-(k-1)\frac{(k+2j)!}{(k+j)!}\right]
=2^{-2(r+1)}\frac{(k+1+2r)!}{(r+1)!(k+r)!}-1
\end{align}

\noindent
By induction we prove equality (\ref{PA Master Equation Sol Final v1  a(n,k) k=2,...,n-2}):\\

\noindent
For $r=0$ we have

\begin{align*}
\frac{2^{-2}}{1!}[2(k-2)-(k-1)]=\frac{2^{-2}(k+1)!}{1!k!}-1
\end{align*}

\noindent
i.e. $2^{-2}(k-3)=2^{-2}(k-3)$.\\

\noindent
By induction hypothesis we prove equality (\ref{PA Master Equation Sol Final v1  a(n,k) k=2,...,n-2}) for $r+1$, supposing it holds for generic $r$, then:

\begin{align*}
&2^{-2(r+1)}\frac{(k+1+2r)!}{(r+1)!(k+r)!}+\frac{2^{-2(r+2)}}{(r+2)!}\left[2(k-2)\frac{(k+1+2r)!}{(k+r)!}-(k-1)\frac{(k+2+2r)!}{(k+1+r)!}\right]
=2^{-2(r+2)}\frac{(k+3+2r)!}{(r+2)!(k+r+1)!}\\&
\\&
2^{2(r+2)}\frac{(k+1+2r)!}{(k+r)!}+2(k-2)\frac{(k+1+2r)!}{(k+r)!}-(k-1)\frac{(k+2+2r)!}{(k+1+r)!}
=\frac{(k+3+2r)!}{(k+r+1)!}
\end{align*}

\newpage
\noindent
i.e.

\begin{align*}
&2(r+2)(k+r+1)+2(k-2)(k+r+1)-(k-1)(k+2+2r)=(k+3+2r)(k+2+2r)
\end{align*}

\noindent
from obvious arrangement w.r.t. $k$ the thesis follow.\\

\noindent
Now consider that $\prod_{i=1}^{n}(2i-1)=\frac{(2n)!}{2^n n!}$, for $k=1,...,n-1$ holds that

\begin{align*}
p_{n,k}&=\frac{(n-1)!}{2^{n-k}\prod_{i=1}^{n}(2i-1)}a_{n,k}=\frac{(n-1)!}{2^{n-k}(2n)!}2^n n! a_{n,k}
\end{align*}

\noindent
i.e.

\begin{align*}
p_{n,k}=\frac{n!(n-1)!}{(2n)!}2^{2n-k}\left[1+(k-1)\sum_{j=0}^{n-k-1}\frac{2^{-2(j+1)}}{(j+1)!}\frac{(k+2j)!}{(k+j)!}\right]
\end{align*}



\end{document}